# Game theory in the case of a geopolitical operation between two persons

O.A. Malafeyev,
Saint-Petersburg State University
malafeyevoa@mail.ru

V. Vekovtsev,
Saint-Petersburg State University

**Abstract**: This article is adjacent to works [1]-[111]. It explores the interaction of two agents during a geopolitical operation. Collaborative work is considered, rather than being done alone. However, each agent has the goal of maximizing personal net profit. We will have 3 different situations depending on the order of the players moves, in each of which we will present the game in both expanded and normal form. In each of them we find optimal strategies for both players and possible Nash equilibrium.

## Section 1. Game conditions

Agent 1 and agent 2 are preparing a geopolitical operation. For the operation to be successful, it is necessary that at least 1 of them personally appear at the operation control center to sign a document with management approving its determination. When this agent 1 is located on another continent, in operations located far from this control center, and therefore the transfer will cost him 4 million dollars. In turn, agent 2 is located in a training center not far from the operation control center, which makes his travel costs negligible. In this case, the potential income from the operation is $20 million. However, the rank of agent 1 is much higher than rank of agent 2, and therefore the main financial actions take place in any case, which is reflected in a larger reward. Both agents can choose whether to plan the course of the operation or go to the operation control center to sign the document.

If both agents come to sign the document, then agent 1 upon completion of the operation will receive $14 million, and agent 2 only $6 million. If agent 1 goes to the operation control center, and agent 2 is planning the operation, then the first will receive already $12 million, and the second $8 million. Finally, if agent 1 remains to plan the operation, and agent 2 goes to the operation control center to sign the document, then they will receive $18 million and $2 million, respectively.

How will agents act if everyone wants to maximize net profit? And an equally important question: who chooses first? We have three situations: 1 – agent 1 chooses first; 2 – agent 2 is selected first; 3 – both agents make a choice at the same time.

## Section 2. Agent 1 chooses first

If we assume that agent 1 will make the first decision, then we get the situation shown in Fig. 1. This tree represents the expanded form of the game. At the top there is a root node (the point labeled "Agent 1") with two branches labeled P (to plan the operation) and S (to go to sign the document). This means that agent 1 has a choice — to take the left branch (P) or the right one (S). This brings us to two nodes with the signature "Agent 2", in each of which agent 2 also has the choice of either planning the operation (P) or going to sign the document (S).

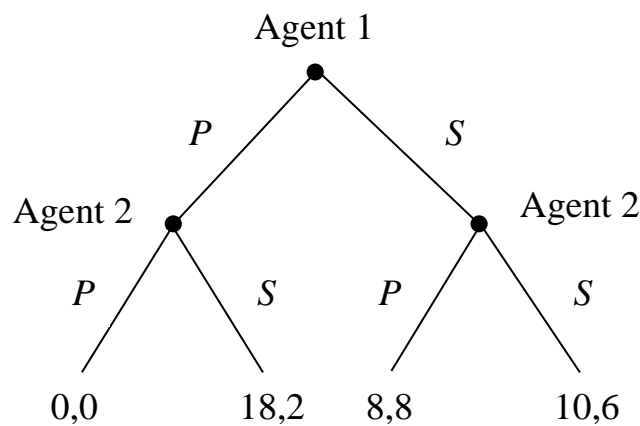

*Figure 1. Expanded form of the game, agent 1 chooses first.*

While agent 1 has only two strategies, agent 2 has four. Let us introduce the corresponding notation:

1. Go to sign the document, regardless of the actions of agent 1 (SS).
2. Engage in planning the operation regardless of the actions of agent 1 (PP).
3. Do the same as agent 1 chooses (PS).
4. Do the opposite of what agent 1 chooses (SP).

Here, the first letter in brackets is responsible for the choice of agent 2 if agent 1 chooses planning, and the second letter is responsible for the choice of agent 1 if agent 1 chooses to go.

Let's call an action a move (choice) made by a player (agent) at a node, and a pure strategy a series of actions that completely determine the player's behavior. Thus, agent 1 has two strategies, each of which is simply his action (P or S). In turn, agent 2 has four strategies, each of which represents two actions: one when agent 2 goes left, and other when he goes right.

At the very bottom of the game tree are four nodes called leaves or terminal nodes. In each of them there is a payoff for both players: the first for agent 1, the second for agent 2, if they both choose strategies that lead to this particular sheet. So, for example, on the leftmost sheet, when both agents choose to schedule an operation, the payoff is (0, 0), since traveling to sign a document with management is a required action for success in the operation. In turn, on the far right sheet, both players choose to go to sign the document. As a result, agent 1 receives $14 million, but his travel costs are $4 million, reducing his net profit to $10 million, and agent 2 receives $6 million, while his costs are negligible. As a result, their net winnings correspond to the specified sheet (10, 6). Likewise with the two central leaves. Thus, the central left sheet implies that agent 1 selects planning, and agent 2 selects travel to sign the document. In this case, agent 1 does not incur costs and receives $18 million, and agent 2, performing an important but easy-to-implement action, ends up with $2 million, which corresponds to the selected sheet (18, 2). Opposite actions are selected by agents in the center right sheet. In this case, agent 1 decides to go to sign the document, incurring costs in the amount of $4 million and receiving 12 million remuneration, which ultimately brings him $8 million in net profit. In turn, agent 2 decides to start planning the operation and, having no costs, receives $8 million, which also corresponds to the sheet in question (8, 8).

On what basis should agent 1 make his choice? Obviously, he needs to find out how agent 2 will react to each of his two options (P or S). If he chooses to schedule the operation, then agent 2 will decide to go to sign the document because he will receive $2 million for this, and not $0 million. Thus, moving along the game tree to the left, agent 1 will receive $18 million. In the case When agent 1 chooses to travel to sign a document, it will be more profitable for agent 2 to choose to plan the operation, receiving $8 million instead of $6 million for this. Thus, moving along the game tree to the right, agent 1 will receive $8 million, as opposed to $18 million in the case movement to the left. Now we know the optimal choice of agent 1 — operation planning.

What can be said in the case of agent 2? Obviously, in the left node he must decide to go to sign the document, but what about the right node? Of course, in this case it doesn't really matter, because he won't be able to end up there anymore. However, it is necessary to indicate not only what the player does "during the game" (in this case, the left branch of the tree), but also all possible nodes of the game tree. This is because we can only say for sure that agent 1 chooses the best answer to agent 2 if we know what choice agent 2 will make, and vice versa. Thus, agent 2 must choose one of the four strategies listed earlier. Obviously, it is profitable for him to play SP, doing the opposite of what agent 1 chooses. This is due to the fact that such a move maximizes his payoff, regardless of agent 1's move.

Conclusion: the only reasonable solution in this game is for agent 1 to plan the operation, and agent 2 to do the opposite of what agent 1 chooses. Their winnings in this case will be $18 million and $2 million, respectively. This outcome is called a Nash equilibrium, which in a two-player game means a pair of strategies, each of which is the best response to the other. In other words, each of these strategies gives the player using it the maximum possible payoff, given the other player's strategy.

There is another way of depicting this game, called strategic or normal form. Typically, both representations or a more convenient one for a particular case are used to solve. The normal form corresponding to Fig. 1 is presented in Fig. 2. In this example, the strategies of the first player (agent 1) are arranged in rows, and the strategies of the second player (agent 2) are arranged in columns. As a result, each entry in the resulting matrix represents the payoffs of both players when they choose the appropriate strategies.

|        | SS   | SP   | PS   | PP   |
|--------|------|------|------|------|
| P      | 18,2 | 18,2 | 0,0  | 0,0  |
| S      | 10,6 | 8,8  | 10,6 | 8,8  |

Agent 2 (columns), Agent 1 (rows: P, S)

*Figure 2. Normal form of the game, agent 1 chooses first.*

In this case, the Nash equilibrium can be found by choosing a row and column such that the payoff at their intersection will be the maximum possible for the first player in the column and the maximum possible for the second in the row. In this case, we can verify that strategy (P, SP) is indeed a Nash equilibrium in a game of normal form, since for agent 1 $18 million is better than $8 million, and $1 million is the best that agent 2 can get in this situation.

Is there another Nash equilibrium in this game? Of course, strategy (P, SS) is likewise a Nash equilibrium because strategy P is the best response to SS and vice versa. But the outcome (P, SS) has the disadvantage that if agent 1 chooses strategy S incorrectly, agent 2 will receive only $6 million, whereas if strategy SP was chosen, he would receive $8 million. In this case, one might say that the strategy SS is weakly dominated by the SP strategy. This means that for agent 2, strategy SP will provide a greater payoff if the opponent chooses strategy S, but if strategy P is chosen, the payoff will remain the same as with SS.

However, what if agent 2, despite everything, does not intend to go to sign the document and chooses the PP or PS strategy? Agent 1 in this case must adhere to strategy S in order to receive $8 million or $10 million instead of $0 million. In turn, agent 2 is more profitable to play PP rather than PS in order to receive $8 million instead of $6 million. As a result, we still have one Nash equilibrium (E, PP), in which agent 2 receives much more ($8 million instead of $2 million), and agent 1 receives much less ($8 million instead of $18 million). This outcome assumes that agent 2 poses an unreliable threat to agent 1. The threat is "inauthentic" because agent 1 knows that if he plays P, then when his opponent's turn comes, he will respond with an SS or SP strategy simply because $2 million is better than $0 million. The Nash equilibrium of a game of expanded form is said to be subgame perfect if its restriction to any subgame of

this game is a Nash equilibrium in it. Thus, strategy (S, PP) is not an ideal subgame, because in the subgame starting with move P of agent 2 on the left in Fig. 1, this is not the best choice.

## Section 3. Junior agent chooses first

Now consider the situation where agent 2 chooses first. The expanded form of this game is shown in Fig. 3. This time agent 1 has 4 strategies (the same ones that belonged to agent 2 in the last game: SS, PP, PS, SP), and agent 2 only has 2 (which belonged to agent 1 in the last game: S and P) .

As in the previous game, to make a decision, agent 2 needs to know how agent 1 will react to both options for his move. If agent 2 chooses to plan an operation, agent 1's best response will be to decide to go to sign the document in order to receive $8 million instead $0 million. And vice versa, if agent 2 chooses a trip, agent 1's best response will be to choose to plan the operation, which will bring him $18 million instead of $10 million. Thus, agent 2 needs to choose strategy P, since in this case he receives $8 million, and when choosing strategy S only $2 million. In this case, agent 1's best choice would be a trip to sign the document, which will ultimately bring both players $8 million. You can note that by choosing first, agent 2 can immediately accept a strategy that posed an unreliable threat when he chose second.

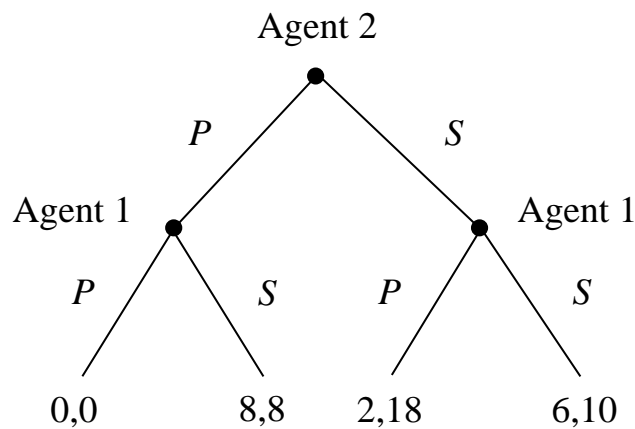

*Figure 3. Expanded form of the game, agent 2 chooses first.*

Let us depict the normal form of the game in the case when agent 2 goes first in Fig. 4. Similarly, we find two Nash equilibria (P, SS) and (P, SP), and also again we discover another one that is not obvious at first glance from the game tree, when already agent 1 poses an unreliable threat to agent 2 with the PP strategy, to which agent 2's best response is move S.

|  |  | Agent 1 |  |  |  |
|---|---|---|---|---|---|
|  |  | *SS* | *SP* | *PS* | *PP* |
| Agent 2 | *P* | 8,8 | 8,8 | 0,0 | 0,0 |
|  | *S* | 6,10 | 2,18 | 6,10 | 2,18 |

*Figure 4. Normal form of the game, agent 2 chooses first.*

# Section 4. Simultaneous selection

The last possible case assumes that both players make a choice at the same time. This is equivalent to a situation where each player makes a decision without knowing about the other player's choice. In this case, everyone has two options: go to the signing of the document (S) or start planning the operation (P). Let's depict this situation in Fig. 5. The dotted line connecting the two nodes that agent 2 selects is called the information set. So, in the general case, an information set is a set of all possible moves for a certain player, taking into account his observations during the game. However, in our case, due to the simultaneous selection of both players, there is the added fact that the players cannot be sure of what has happened at this point in the game and what their position is. Thus, you can safely swap agent 1 and agent 2, of course, changing their winnings at the terminal nodes. Ultimately, this means that there can be more than one extended form game representing the same game situation.

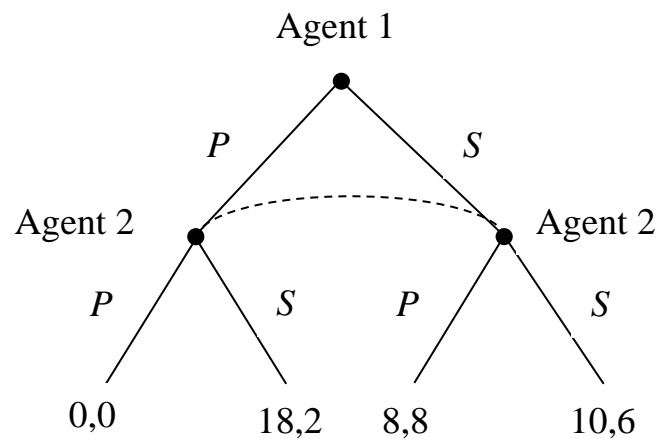

*Figure 5. Expanded form of the game, simultaneous choice.*

Even though there is less strategy in this game, looking at the game tree it is difficult to see what the equilibrium might be. This happens because the choice of agent 2 cannot depend on what choice agent 1 makes, since agent 2 has no information about this. Consider the normal form game presented in Fig. 6. It is easy to see that both strategies, (P, S) and (S, P), are Nash equilibria, the first of which favors agent 1, and the second favors agent 2. Thus, in the first case, the agents receive $18 million and $2 million respectively, and in the second $8 million each. There is also a not so obvious Nash equilibrium in the case where both agent 1 and agent 2 randomly choose P or S with probability 1/2. This is called a mixed strategy Nash equilibrium. In this case, agent 1 will receive $9 million, and agent 2 will be content with $4 million. The gains are not so great due to the fact that with probability 1/4 both agents choose to plan the operation and receive zero profit.

|  | Agent 2 | |
|---|---|---|
|  | S | P |
| Agent 1  S | 10,6 | 8,8 |
| P | 18,2 | 0,0 |

*Figure 6. Normal game form, simultaneous choice.*

# Conclusion

In the article is formalized and studied the problem of maximizing personal profit in a game of two agents. Three different cases are considered depending on the order of choice. In each situation, Nash equilibria and optimal strategies for both agents are found. The corresponding payoffs are found and both expanded and normal forms of the game are constructed.